\documentclass[12pt]{amsart}

\usepackage{amssymb,array,epsfig}
\usepackage{color}
\newtheorem{lemma}{Lemma}

\newtheorem{theorem}[lemma]{Theorem}

\newtheorem{remark}[lemma]{Remark}

\newcommand{\R}{{\mathbb R}}
\newcommand{\C}{{\mathbb C}}

\newcommand{\Dom}{{\rm Dom}}

\newcommand{\Spec}{{\rm Spec}}

\renewcommand{\Re}{{\rm Re}\;}
\renewcommand{\Im}{{\rm Im}\;}

\newcommand{\dist}{{\rm dist}}

\newcommand{\dsc}{\mathrm{disc}}
\newcommand{\ess}{\mathrm{ess}}

\title[Stability of Quadratic Projection Methods]
{Stability of Quadratic Projection Methods}

\author[L.~Boulton and M.~Strauss]{Lyonell Boulton \and
Michael Strauss}

\date{25th January 2007}



\begin{document}

\begin{abstract}

We discuss stability properties of the method studied recently in
\cite{lesh} and \cite{bo2},
for computing eigenvalues in gaps of the essential spectrum.

\end{abstract}

\maketitle

\newpage


\section{Introduction}
\label{sec0}

\subsection{Spectral Pollution in the Galerkin method}
\label{ss1.1}
Let $A$ be a self-adjoint operator acting on an infinite dimensional
 Hilbert space
$\mathcal{H}$, with a dense domain $\Dom (A)$. The spectrum of $A$,
$\mathrm{Spec}(A)$, may be expressed as the union of the discrete
spectrum consisting of all isolated eigenvalues of finite
multiplicity, $\mathrm{Spec}_{\dsc}(A)$, and the essential spectrum,
where
$\mathrm{Spec}_{\ess}(A):=\mathrm{Spec}(A)\backslash\mathrm{Spec}_{\dsc}(A)$.
In most standard situations the essential spectrum can be found
analytically, but points in $\Spec_{\dsc}(A)$ are usually estimated
by numerical procedures.

The estimation of $\Spec_{\dsc}(A)$ is often performed through
subspaces $\mathcal{L}\subset \Dom (A)$ and
corresponding truncations of $A$. Standard numerical techniques, such as the finite element method, aim at
solving Galerkin approximate problems posed
in weak form:

\medskip
\ \hspace{.3in} (P) \hspace{1cm}
\begin{minipage}[c]{3in}
find $0\not=u\!\in\!\mathcal{L}$ and
$\lambda\!\in\!\mathbb{R}\,$ such that \newline $\langle Au,v\rangle
=\lambda\langle u,v\rangle \qquad \forall v\in \mathcal{L}$
\end{minipage}
\medskip

\noindent where $\mathcal{L}$ is finite dimensional.

Backed by the Rayleigh-Ritz variational principle, when applicable,
the Galerkin method represents a powerful tool in the analysis of
spectral properties of linear operators. However, the Galerkin
method is not foolproof, in general, the solutions of (P) might fail
to provide reliable information about the location
of $\mathrm{Spec}(A)$ (see \cite{ds}, \cite{da0}, \cite{lesh},
\cite{dapl}, \cite{sp2}).

The drawbacks in the Galerkin method are due in part to the so
called spectral pollution phenomenon which we now describe. Let
$\mathcal{L}_n\subset \Dom (A)$ be a sequence of subspaces
approaching
$\mathcal{H}$, as $n\to\infty$  (e.g.
satisfying \eqref{g1} below with $p=0,1$ only).
Suppose we found $0\not=u_n\in \mathcal{L}_n$ and $\lambda_n\in
\mathbb{R}$ solutions of (P) with $\mathcal{L}=\mathcal{L}_n$,
satisfying $\lambda_n\to \mu$ and $\|u_n\|^{-1}u_n\to w \in \Dom(A)$
in the weak topology. By the approximating property of
$\mathcal{L}_n$, we may obtain
\[
    \langle Aw-\mu w,v\rangle
=0
\qquad \forall v\in \mathrm{Dom}(A),
\]
which appears to suggest that $\mu$ is in $\mathrm{Spec}(A)$.
Unfortunately, the latter conclusion is not ensured in general. Without
further information about the structure of $A$ (e.g. compactness
properties), $w$ might be $0\in \Dom (A)$, so convergent solutions
of the approximate problem might produce ``polluted'' sequences
$\lambda_n\to\mu\not \in \mathrm{Spec}(A)$.

The emergence of spurious eigenvalues in gaps of
$\Spec_{\ess}(A)$ represents a serious difficulty in applications
such as elasticity theory and solid state physics (see
\cite{ds} and \cite{sp2}), as there is no universal
recipe to detect or prevent them for a given operator
$A$ and sequence of approximate subspaces $\mathcal{L}_n$.


\subsection{Pollution-free strategies and quadratic methods}
\label{ss1.2}
Spectral pollution is a consequence of the fact
that in (P) we are truncating simultaneously both $u$ and $v$. Indeed, let
$\Pi$ be the orthogonal projection onto $\mathcal{L}$ and
\[
\hat{F}_{\mathcal{L}}(x):= \min_{0\not= v\in \mathcal{L}}
 \frac{\|\Pi(x-A)v\|}{\|v\|}.
\]
Then $\hat{\lambda} \in \mathbb{R}$ satisfies (P) if, and only if,
$\hat{F}_{\mathcal{L}}(\hat{\lambda})=0$. That is to say, there
exists $\hat{u}\in \mathcal{L}$ such that
$(\hat{\lambda}-A)\hat{u}\perp \mathcal{L}$. As
$\|(\hat{\lambda}-A)\hat{u}\|/\|\hat{u}\|$ is not guaranteed to be
small, we have no indication whether $\hat{\lambda}$ is close to
$\mathrm{Spec}(A)$ or not.

This argument suggests that the correct
quantity to look at is
\[
     F_{\mathcal{L}}(x):=
    \min_{0\not= v\in \mathcal{L}} \frac{\|(x-A)v\|}{\|v\|}.
\]
As
\[
    F_{\mathcal{L}}(x)\geq \inf_{u\in \Dom (A) }
 \frac{\|(x-A)u\|}{\|u\|} = \|(x-A)^{-1}\|^{-1}=
 \mathrm{dist}\,[x,\mathrm{Spec}(A)],
\]
$F_{\mathcal{L}}(x)$ can be
 close to $0$ only when $x$ is close to a point in
the spectrum of $A$.

In \cite{dapl}, Davies and Plum considered a pollution-free strategy
for finding $\Spec(A)$ based on computing the profile of $F_{\mathcal{L}}(x)$
for $x\in \mathbb{R}$. If $\mathcal{L} \subset \Dom (A^2)$,
\begin{equation} \label{e13}
\begin{aligned}
 F_{\mathcal{L}}(x)^2 & = \min_{0\not= v\in \mathcal{L}}
     \frac{\langle \Pi(x-A)^2 v,v\rangle }{\|v\|^2}
        =
 \|[\Pi(x-A)^2\upharpoonright \mathcal{L}]^{-1}\|^{-1} \\
 & = \min_{0\not= v\in \mathcal{L}}
 \frac{\|\Pi(x-A)^2 v\|}{\|v\|}
 =: G_{\mathcal{L}}(x).
\end{aligned}
\end{equation}
Therefore estimating $F_{\mathcal{L}}(x)$ reduces to computing eigenvalues of
self-adjoint matrices depending on the parameter $x\in\mathbb{R}$.

The approach developed in \cite{dapl} relies heavily on being able
to find accurately a matrix representation for
$\Pi(x-A)^2\upharpoonright\mathcal{L}$ in terms of an orthonormal
basis of $\mathcal{L}$. This is a drawback, for instance, if
$\mathcal{L}$ is given by the finite element method, where an
orthonormalisation of the basis will be numerically expensive.

An alternative pollution-free method which is independent of the
matrix representation of $\Pi(x-A)^2\upharpoonright\mathcal{L}$ is
also available and it may be obtained by considering the zeros of
the function $G_{\mathcal{L}}(z)$ for $z\in \C$. Typically
$G_{\mathcal{L}}(z)$ and $F_{\mathcal{L}}^2(z)$ only coincide at
$z\in \R$. The $(2\dim \mathcal{L})$ zeros of the polynomial
$\det\big(\Pi(z-A)^2\upharpoonright \mathcal{L}\big)$ are the zeros
of $G_{\mathcal{L}}(z)$ and, on the other hand,
$F_{\mathcal{L}}(z)\not=0$ unless $z$ is an eigenvalue of $A$, with
corresponding eigenvector $u\in\mathcal{L}$, a very unlikely
situation. It is remarkable, however, that the non-real zeros of
$G_{\mathcal{L}}(z)$ also provide reliable information about the
location of $\Spec (A)$.

This alternative procedure has been recently discussed in
\cite{lesh}, \cite{bo1} and
\cite{bo2}, and it can be traced back to \cite{da0} and \cite{shar}.
A central role is played by the problem

\medskip

\ \hspace{.3in}(Q) \hspace{.2in}\begin{minipage}[c]{4in}find
$\zeta\!\in\! \mathbb{C}$ such that $\exists u\in \mathcal{L}$ with
\newline $\langle Au,Av\rangle-2\zeta\langle Au,v\rangle+
\zeta^2 \langle u,v\rangle =0, \quad \forall v\in \mathcal{L}.$
\end{minipage}

\medskip

\noindent It is readily seen that $G_{\mathcal{L}}(\zeta)=0$  if, and only if, $\zeta$
is a solution of (Q). The philosophy of the method is to regard (Q),
in place of (P), as an approximate spectral problem for operator
$A$.

The following universal non-pollution result justifies favouring
(Q) over (P) (see \cite[Theorem~2.6]{lesh} or Theorem~\ref{pract}
below): if $\zeta$ is a solution of (Q), then
\begin{equation} \label{eu}
\mathrm{dist}[\Re \zeta,\mathrm{Spec}(A)]\leq |\Im \zeta|.
\end{equation}
That is to say, $\zeta$ can be close to $\R$, \emph{only when} it is
also close to the spectrum of $A$.

Problem (Q) gives rise to a
matrix spectral problem quadratic in the spectral parameter. This
added complication balances out with the reliability of the method
expressed in the above result.

Now, will a solution of (Q) \emph{ever} be close to $\R$? As for the
Galerkin method, in general, additional conditions on a sequence of
subspaces $\mathcal{L}_n$ are required for convergence. A precise
statement reads as follows, see \cite{bo2} or Theorem~\ref{t3} below.
Let $\lambda\in
\Spec_{\dsc}(A)$ and $\Pi_n$ be the
orthogonal projection onto $\mathcal{L}_n\subset \Dom (A^2)$. If
\begin{equation} \label{g1}
\frac{\|\Pi_n A^p\Pi_nu-\lambda^pu\|}{\|u\|}\to 0,
\qquad
 \begin{array}{ll} \forall p=0,1,2, \\
 \forall u\in \Dom(A): Au=\lambda u, \end{array}
\end{equation}
then there exists $\zeta_n\in\C$
satisfying (Q) with $\zeta=\zeta_n$ and $\mathcal{L}=\mathcal{L}_n$,
such that $|\zeta_n-\lambda|\to 0.$ The above hypothesis is fulfilled
immediately, for instance, if $A$ is bounded and $\Pi_n v\to v$ for
all $v\in \mathcal{H}$.

The combination of these two results appears to provide a general
pollution-free procedure for finding discrete eigenvalues of
self-adjoint operators. Although this might seem too optimistic at
the present moment, one of the advantages of this method lies in the
fact that it is applicable without any special restriction upon the
structure of $\Spec (A)$. Moreover, the  requirements on
$\mathcal{L}_n$  are analogous to those needed in the Galerkin
method.


\subsection{Stability of Quadratic Projection Methods}
On the downside, here we are confronted with a more difficult
problem to solve. In general, the finite-dimensional eigenvalue
problem associated to (Q) is non-Hermitian. Accuracy, as well as
stability of the method becomes a delicate matter. The main goal of
the present note is to discuss how non-pollution and
convergence of the method are affected, when the coefficients of
problem (Q) are known only approximately.

In Section~\ref{sec2} we will show that the non-pollution  property
remains stable in a sense which will be specified below. In
Section~\ref{sec3} we will discuss stability of approximation.
 Note that a consistent formulation of (Q) only requires
$\mathcal{L}_n\subset \Dom (A)$, see Remark~\ref{r1}. Under a suitable hypothesis on the
subspaces $\mathcal{L}_n$, our Theorem~\ref{t3} extends the
analogous result of \cite{bo2} by allowing $\mathcal{L}_n \cap [\Dom
(A) \setminus \Dom (A^2)]\not= \varnothing$. In the final section
we report on various numerical experiments performed on a simple
example.


\section{Pollution-free Stability}
\label{sec2}

We devote this section to showing that, given error bounds in the
computation of the coefficients of problem (Q), it is possible to
control errors in the pollution-free estimation of $\Spec (A)$ by
the quadratic method described in
Section~\ref{ss1.2}.

Let us begin by fixing some notation.
Below $\mathcal{L}_n$ is an $n$-dimensional subspace of
$\mathrm{Dom}(A)$ with basis $\{e_1,\dots,e_n\}$. This basis
will always be normalised, $\|e_j\|=1$ for all $j=1,\ldots,n$.
When sufficiently clear from the context, we will suppress
the sub-index and write $\mathcal{L}\equiv \mathcal{L}_n$.

For any $u\in\mathcal{L}$,
$u = \overline{u_1}e_1 +\dots+\overline{u_n}e_n,$
from which we define the following norm on
$\mathcal{L}$,
\begin{equation*}
\Vert u\Vert_0 := \big{(}\vert u _1\vert^2 + \dots
+\vert u_n\vert^2\big{)}^{\frac{1}{2}}.
\end{equation*}
Since $\mathcal{L}$ is a finite dimensional space, there exists
$\beta > 0$, such that
\begin{equation} \label{e22}
\Vert u \Vert = \langle u , u \rangle^{\frac{1}{2}} \ge \beta\Vert u
\Vert_0\quad\forall u \in\mathcal{L}.
\end{equation}
If $\{e_1,\ldots,e_n\}$ is an orthonormal basis, then
$\|\!\cdot\!\|= \|\!\cdot\! \|_0$. However when the basis is far
from being orthonormal, $\beta$ will be small. We will occasionally
write $\mathbf{u} = (u_1,\ldots,u_n)\in \mathbb{C}^n$.

Let matrices $A_0$, $A_1$ and $A_2$ in $\mathbb{C}^{n\times n}$
be given entrywise by
\begin{equation} \label{e10}
   [A_0]_{jk}=\langle A e_j,A e_k \rangle ,\quad
   [A_1]_{jk}=\langle A e_j,e_k \rangle ,
   \quad [A_2]_{jk}=\langle e_j,e_k \rangle .
\end{equation}
Define the matrix polynomial
$M(z)\in \mathbb{C}^{n\times n}$ as
\begin{align} \label{e21}
M(z): = A_0 - 2zA_1 + z^2A_2, \qquad \qquad z\in\mathbb{C}.
\end{align}
Then $\zeta\in\C$ is a solution of (Q) if, and
only if, $\det [M(\zeta)]=0$.

The stability results we establish below give a
positive answer to the following question. Suppose we are only able
to estimate the matrices $A_p$ by $\tilde{A}_p$ and the norm of the
error $\|A_p-\tilde{A}_p\|\leq \varepsilon_p$, $p=0,1,2$. Can we
recover information about the spectrum of $A$ from the approximate
problem

\smallskip
\ \hspace{.3in} ($\tilde{\mathrm{Q}}$) \hspace{.2in}
\begin{minipage}[c]{3.5in}
find $\zeta\!\in\! \mathbb{C}\,:$ $\det [ \tilde{A}_0-2\zeta
\tilde{A}_1 + \zeta^2\tilde{A}_2]=0$
\end{minipage}

\smallskip

\noindent with accuracy possibly
depending upon $\varepsilon_p$?

\medskip

The following preliminary result will be needed.

\begin{lemma}
For $z\in\mathbb{C}$ and $\delta > 0$,
let
\begin{eqnarray*}
J &=& [\Re{z} - \vert\Im{z}\vert - \delta ,\Re{z} + \vert\Im{z}\vert + \delta ],\\
\Omega &=& \big{\{}(x - z)^2 : x\in\mathbb{R}\backslash J\big{\}}.
\end{eqnarray*}
Then,
\begin{equation}\label{1}
\inf_{\upsilon\in \Omega}\Re{\upsilon} = 2\delta\vert\Im{z}\vert +
\delta^2 .
\end{equation}
\end{lemma}
\begin{proof}
Let $z = a + ib$, $a,b\in \R$, then for any $\upsilon\in \Omega$ we have
\begin{displaymath}
\upsilon = (x - a)^2 - b^2 - 2b(x - a)i
\end{displaymath}
for some $x\in\mathbb{R}\backslash J$. It is clear
that $\Re{\upsilon} > 0$ and moreover
\begin{eqnarray*}
\inf_{\upsilon\in \Omega}\Re{\upsilon} &=& \inf_{x\in\mathbb{R}\backslash J} (x - a)^2 - b^2\\
&=& (\vert b\vert + \delta)^2 - b^2\\
&=& 2\delta\vert b\vert + \delta^2
\end{eqnarray*}
verifying \eqref{1}.
\end{proof}

\begin{theorem}\label{pert}
Let $A$ be a self-adjoint operator acting on a Hilbert space
$\mathcal{H}$, and $\mathcal{L}$ be an $n$-dimensional subspace of
$\Dom (A)$. Let $B$ be a singular $n\times n$ matrix. For any
$z\in\mathbb{C}$, let $M(z)$ and $\beta$ be as in \eqref{e21} and
\eqref{e22}. Let $\alpha_{z}\in\mathbb{R}$ with $\alpha_{z}\ge\Vert
M(z) - B\Vert_{\mathbb{C}^n}$. If $\delta>0$ is such that
$(\delta^2+2\delta |\Im{z}|)\beta^2>\alpha_z$, then
 \begin{equation}\label{2}
 \Spec(A)\cap [\Re{z} - \vert\Im{z}\vert - \delta ,\Re{z} +
\vert\Im{z}\vert + \delta ]\ne\varnothing.
 \end{equation}
\end{theorem}
\begin{proof}
Let $\delta>0$ be as in the hypothesis and suppose the intersection
\eqref{2} is empty. Using the spectral theorem and \eqref{1}, we
have for all $u\in\mathcal{L}$
\begin{align*}
\Re  & (\overline{\mathbf{u}}^\mathrm{T} M(z) \mathbf{u})
 = \Re\sum_{jk=1}^n\langle
(A - z)e_j,(A - \overline{z})e_k\rangle u _k\overline{ u }_j \\ &=
\Re\sum_{jk=1}^n\langle (A - z)\overline{ u }_je_j,(A -
\overline{z})\overline{ u }_ke_k\rangle = \Re\langle(A - z) u ,(A -
\overline{z}) u \rangle\\ &= \int_{\mathbb{R}}\Re(\lambda -
z)^2\textrm{ }d\langle E_\lambda u,u\rangle \ge
(2\delta\vert\Im{z}\vert + \delta^2)\Vert u \Vert^2 \\
&\ge (2\delta\vert\Im{z}\vert + \delta^2)\beta^2\Vert u  \Vert_0^2 =
(2\delta\vert\Im{z}\vert + \delta^2)\beta^2\Vert \mathbf{u}
\Vert_{\mathbb{C}^n}^2,
\end{align*}
where $E_\lambda$ is the spectral measure associated to $A$.
It then follows from the Schwarz inequality that for any $ \mathbf{u}
\in\mathbb{C}^n$
\begin{eqnarray*}
\Vert M(z) \mathbf{u} \Vert_{\mathbb{C}^n} &\ge&
(2\delta\vert\Im{z}\vert + \delta^2)\beta^2\Vert \mathbf{u}
\Vert_{\mathbb{C}^n},
\end{eqnarray*}
so that the operator $M(z) : \mathbb{C}^n\rightarrow\mathbb{C}^n$ is
invertible and
\begin{eqnarray*}
\Vert M(z)^{-1}\Vert_{\mathbb{C}^n} &\le& \big{(}(2\delta\vert
\Im{z}\vert + \delta^2)\beta^2\big{)}^{-1} < \alpha_{z}^{-1}.
\end{eqnarray*}
In particular $\Vert M(z)^{-1}\Vert_{\mathbb{C}^n}^{-1} > \Vert M(z)
- B\Vert_{\mathbb{C}^n}$, from which it follows that $B$ is
not singular. The result follows from the obtained contradiction.
\end{proof}

The next theorem is the main result of this section and it is
an improvement on \cite[Theorem~2.6]{lesh}.

\begin{theorem}\label{pract}
Let $A$ be a self-adjoint operator acting on a Hilbert space
$\mathcal{H}$, and $\mathcal{L}$ be an $n$-dimensional subspace of
$\Dom (A)$. Define $A_0,\,A_1$ and $A_2$ as in \eqref{e10}. Let
$\tilde{A}_{p}$ be $n\times n$
matrices, such that for
$\varepsilon_{p}\ge 0$
\begin{eqnarray*}
\Vert A_p - \tilde{A}_{p}\Vert_{\mathbb{C}^n} \leq \varepsilon_{p},
\qquad p=0,1,2.
\end{eqnarray*}
If the matrix $\tilde{A}_{0} - 2\zeta \tilde{A}_{1} + \zeta^2\tilde{A}_{2}$
is singular for some $\zeta\in\mathbb{C}$, then
\begin{equation}\label{6}
\Spec(A)\cap [\Re{\zeta} - \tilde\delta ,\Re{\zeta} +
\tilde\delta]\ne\varnothing
\end{equation}
for
\begin{equation*}
\tilde\delta = \sqrt{\vert\Im{\zeta}\vert^2 + \beta^{-2}(\vert
\zeta\vert^2\varepsilon_2 + 2\vert \zeta\vert\varepsilon_{1} +
\varepsilon_{0})}.
\end{equation*}
\end{theorem}
\begin{proof}
With the notation of Theorem \ref{pert}, take $B = \tilde{A}_{0}
-2\zeta\tilde{A}_{1} + \zeta^2\tilde{A}_{2}$. Since
\begin{displaymath}
\Vert M(\zeta) - B\Vert_{\mathbb{C}^n} \leq (\vert \zeta\vert^2\varepsilon_2
+ 2\vert \zeta\vert\varepsilon_{1} + \varepsilon_{0}),
\end{displaymath}
\eqref{6} follows from Theorem \ref{pert}.
\end{proof}

In particular, under the hypothesis above,
\begin{equation}\label{6a}
\mathrm{dist}[\Re \zeta,\Spec (A)]\leq \vert\Im \zeta\vert +\beta^{-1}\epsilon
\end{equation}
with $\epsilon = \sqrt{\vert
\zeta\vert^2\varepsilon_2 + 2\vert \zeta\vert\varepsilon_{1} +
\varepsilon_{0}}$.
If the basis $\{e_1,\dots,e_n\}$ is orthonormal, then \eqref{6} and
\eqref{6a} hold with $\beta =1$. Note that the case $\varepsilon_0=
\varepsilon_1=\varepsilon_2=0$ corresponds to
\cite[Theorem~2.6]{lesh}, see \eqref{eu}.

\section{Stability of Convergence in the Quadratic Method}
\label{sec3}

A consistent formulation of (Q) only requires $\mathcal{L}\subset
\Dom(A)$. However, the available approximation results for the
quadratic method (cf. \cite{bo1} and \cite{bo2}) impose the
hypothesis $\mathcal{L}\subset \Dom(A^2)$.
 In this
section we show that, if  $\mathcal{L}_n\subset \Dom(A)$ approach
reasonably well the eigenspace associated to an eigenvalue
$\lambda\in \Spec_{\dsc}(A)$, then solutions of (Q) will
converge to $\lambda$ in the large $n$ limit, and the process
remains stable under perturbation of the matrix coefficients
of the polynomial $M(z)$.

\begin{remark} \label{r1}
Allowing the possibility of test spaces $\mathcal{L}\nsubseteq
Dom(A^2)$ is only relevant when $A$ is unbounded. If $A$ is a
differential operator of order $2m$ and the trial spaces are
constructed using the finite element method, $\mathcal{L}\subset
\Dom(A^2)$ requires $C^{4m-1}$ conforming elements, while
$\mathcal{L}\subset \Dom(A)$ only requires $C^{2m-1}$ conforming
elements. The performance of the interpolation algorithm in the
finite element method is usually compromised as $m$ increases.
\end{remark}

Below we highlight explicitly the dependency on $n$ of approximate
subspaces and operators, so we denote matrices $M(z)$ and $A_p$,
corresponding to $\mathcal{L}_n$, by $M^{(n)}(z)$ and $A^{(n)}_p$,
respectively. We also assume throughout this section that the basis
$\{e_1,\ldots,e_n\}$ of $\mathcal{L}_n$ is orthonormal. In general
we do not assume that $\mathcal{L}_n\subseteq\mathcal{L}_m$ whenever
$n<m$. Strictly speaking we should denote the basis functions of
$\mathcal{L}_n$ by $\{e_j^{(n)}\}$. However we suppress this
notation as no confusion shall arise.

For $u\in \Dom (A)$, the projection of $u$ onto $\mathcal{L}_n$ is
then given by
\[
      \Pi_n u=\sum _{k=1}^n \langle u,e_k\rangle e_k=
     \sum _{k=1}^n \overline{u}_k e_k.
\]
Since $\{e_1,\ldots,
e_n\}$ is  orthonormal, we can isometrically identify $\mathcal{L}_n$ with
$\C^n$.

Our key result assumes the following hypothesis on the sequence
$\mathcal{L}_n$:
\begin{equation*}(\mathrm{H}) \qquad  \begin{array}{l}
\forall p,q=0,1, \textrm{ and } \forall u\in \Dom (A)
: Au=\lambda u,   \\
 \|\sum_{j=1}^n \langle A^p \Pi_n u,A^q e_j \rangle e_j -\lambda^{p+q}u
\|\to 0,\textrm{ as }n\to\infty.
\end{array}
\end{equation*}
Whenever $\mathcal{L}_n\subset \Dom(A^2)$, (H)
reduces to \eqref{g1}. Furthermore, if $A$ is bounded and
$\Pi_n$ converges strongly to the identity, then (H)
holds true for all $\lambda\in \Spec_{\dsc}(A)$.

The following result is an improvement upon \cite[Theorem~2.2]{bo2}.

\begin{theorem} \label{t3}
Let $A$ be a self-adjoint operator on a Hilbert space.
Suppose that the sequence of approximate subspaces
$\mathcal{L}_n\subset \Dom (A)$ satisfy (H). Let $\lambda\in
\Spec_{\dsc} (A)$ and let $d:=\dist [\lambda,
\Spec(A)\setminus \{\lambda \}]$. Given $0<\delta< d/4$, there
always exist
$ N ,\,\varepsilon_0,\varepsilon_1, \varepsilon_2>0$ ensuring
the following. If $n> N $ and the matrices $\tilde{A}_p \in
\C^{n\times n}$ satisfy
\[
    \|\tilde{A}_p-A^{(n)}_p \| <\varepsilon _p,
\qquad p=0,1,2,
\]
then
\begin{itemize}
\item[(a)] we can always find $\zeta\in \C$ with $
    \det [\tilde{A}_0-2\zeta\tilde{A}_1
   +\zeta^2\tilde{A}_2]=0$
and $|\zeta-\lambda |<\delta$,
\item[(b)] the set $\{\mu \in \C \,:\,
\det [\tilde{A}_0-2\mu\tilde{A}_1
   +\mu^2\tilde{A}_2]=0 \}$ does not intersect
the annulus $\{w\in\C\,:\,\delta < |w-\lambda|\leq d/4 \}$.
\end{itemize}
\end{theorem}

The proof of this result will be given at the end
of this section. It will be a consequence of various
technical lemmas, in particular, suitable extensions
of \cite[Lemmas~5.1 and 5.3]{bo2} and
various regularity properties of $G_{\mathcal{L}}(z)$.

\medskip

We begin with the rigorous definition of the right hand side
of \eqref{e13} in the case $\mathcal{L}\subset \Dom(A)$.
For $z\in \C$, let
\begin{equation*}
G_{\mathcal{L}}(z):= \min_{0\not = \mathbf{v} \in\C^n}
\frac{\|M(z)\mathbf{v}
\|_{\mathbb{C}^n}}{\|\mathbf{v}\|_{\mathbb{C}^n}}.
\end{equation*}
If $\mathcal{L}\subset \Dom(A^2)$, then $G_{\mathcal{L}}(z)$ coincides with
the right hand side of \eqref{e13}. Below we will write
$G_n(z)\equiv G_{\mathcal{L}_n}(z)$.

Clearly $G_{\mathcal{L}}(\zeta)=0$
if, and only if, $\det M(\zeta)=0$; so the solutions of problem (Q)
are completely characterised as the zeros of $G_{\mathcal{L}}(z)$.
It is readily seen that:
\begin{equation} \label{e15}
  G_{\mathcal{L}}(z)= \| [M(z)]^{-1} \|^{-1} =
  \mathrm{least\ singular\ value\ of}\,M(z).
\end{equation}
In fact $G_{\mathcal{L}}(z)^{-1}$ is a continuous subharmonic
function in the region $\{z\in \C:\det
M(z)\not=0\}$, with singularities at the zeros of
$\det[M(z)]$ (see e.g. \cite{da0} or \cite[Lemma~4.1]{bo2}). This
property will play a central role below.

\medskip

The statement of Theorem~\ref{t3} will be obtained as a consequence
of the fact that $G_{\mathcal{L}}(z)$ is small if,
and only if, for small enough $\varepsilon_p>0$, $z=\zeta$ is
a solution of an approximate problem ($\tilde{\mathrm{Q}}$).
The following notion, which has recently become standard,
will simplify considerably most of our arguments. Let
\[
  \Lambda^{\mathcal{L}}(\varepsilon_0,\varepsilon_1,\varepsilon_2):=
\{ z\in \C\,:\, G_{\mathcal{L}}(z) -(\varepsilon_0+2 \varepsilon_1 |z|+
   \varepsilon_2 |z|^2) \leq 0 \}.
\]
This set is called
the structured pseudospectrum of the matrix polynomial
$M(z)$, see \cite{hiti}.

The proof of the following fundamental property of the
pseudospectrum is a direct consequence of
\eqref{e15} and \cite[Lemma~2.1]{hiti}. It clearly suggests how to
verify the validity of (a) and (b) of Theorem~\ref{t3}.

\begin{lemma} \label{e3}
The complex number $\zeta\in \Lambda^{\mathcal{L}}(\varepsilon_0,\varepsilon_1,\varepsilon_2)$ if, and only if, $\det [\tilde{A}_0-2\zeta\tilde{A}_1
   +\zeta^2\tilde{A}_2]=0$ for some $\tilde{A}_p\in \C^{n\times n}$
satisfying  $\|\tilde{A}_p-A_p \| \leq \varepsilon_p, \ p=0,1,2.$
\end{lemma}

Furthermore, cf. \cite[Theorem~2.3]{laps},

\begin{lemma} \label{t5}
Let $\Omega$ be a connected
component of
$\Lambda^{\mathcal{L}}(\varepsilon_0,\varepsilon_1,\varepsilon_2)$, such that
$\det M(\mu)=0$ for some $\mu \in \Omega$. If
$\|\tilde{A}_p-A_p \| \leq \varepsilon_p$, there always exist
$\zeta \in \Omega$ such that $\det [\tilde{A}_0-2\zeta\tilde{A}_1
   +\zeta^2\tilde{A}_2]=0$.
\end{lemma}

\medskip

We now establish two key relations between
the large $n$ limit of $G_n(z)$ and $\mathrm{dist} [z,\Spec(A)]^2$
in a neighbourhood  of the discrete spectrum of $A$.

\begin{lemma} \label{t1} Let $\lambda\in \Spec_{\dsc} A$.
If the sequence of approximate subspaces
$\mathcal{L}_n\subset \Dom (A)$ satisfy (H), then
\[
\lim_{n\to \infty} G_n(\lambda) =0.
\]
\end{lemma}
\begin{proof}
Let $u\in \Dom (A)\backslash\{0\}$ be such that
$Au=\lambda u$. Consider the vector $\mathbf{u} =
(\langle e_1,u\rangle,\dots,\langle e_n,u\rangle)$. We have
\begin{eqnarray*}
[M^{(n)}(\lambda)\mathbf{u}]_i &=& \sum_{j=1}^n[M(\lambda)]_{ij}\langle e_j,u\rangle\\
&=& \sum_{j=1}^n\langle Ae_i,A\langle u,e_j\rangle e_j\rangle -
2\lambda\langle Ae_i,\langle u,e_j\rangle e_j\rangle +
\lambda^2\langle e_i,\langle u,e_j\rangle e_j\rangle\\
&=& \langle Ae_i,A\Pi_nu\rangle - 2\lambda\langle Ae_i,\Pi_nu\rangle
+ \lambda^2\langle e_i,\Pi_nu\rangle,
\end{eqnarray*}
so that
\begin{eqnarray*} \Vert
M^{(n)}(\lambda)\mathbf{u}\Vert_{\mathbb{C}^n}^2 &=&
\sum_{j=1}^n\vert\langle A\Pi_nu,Ae_j\rangle - 2\lambda\langle
\Pi_nu,Ae_j\rangle +
\lambda^2\langle \Pi_nu,e_j\rangle\vert^2\\
&=& \Vert\sum_{j=1}^n\langle A\Pi_nu,Ae_j\rangle e_j - 2\lambda\langle
\Pi_nu,Ae_j\rangle e_j + \lambda^2\langle \Pi_nu,e_j\rangle
e_j\Vert^2.
\end{eqnarray*}
The right hand side converges to zero
by virtue of (H).
Also $\Vert\mathbf{u}\Vert_{\mathbb{C}^n}\to\Vert u\Vert\ne 0$ as
$n\to\infty$ . Now, fix $\varepsilon>0$. Then, for
all $n$ large enough,
\begin{align*}
   G_n(\lambda) & \leq \frac{\|M^{(n)}(\lambda)\mathbf{u}\|_{\mathbb{C}^n}}
   {\|\mathbf{u}\|_{\mathbb{C}^n}} \leq \varepsilon.
\end{align*}
As $G_n(z)$ is non-negative and $\varepsilon$ is arbitrary
the lemma follows.
\end{proof}

In general it is possible to construct examples where $\lim_{n\to
\infty} G_n(z)=0$ for certain
$z\not\in \R$, \cite{bo1}. However, this is not possible
for $z$ in the vicinity of discrete eigenvalues of $A$.

\begin{lemma} \label{t2}
Let $\lambda\in \Spec_{\dsc} (A)$ and let $d>0$ be as in
Theorem~\ref{t3}. Assume that the sequence of approximate subspaces
$\mathcal{L}_n\subset \Dom (A)$ satisfy (H). For all
$0<\delta<d/4$, there exist a constant $0<s\leq 1$
such that
\begin{equation} \label{e14}
    \liminf_{n\to \infty} G_n(z) \geq s \delta^2 \qquad \mathrm{for\ all} \qquad
    \delta\leq |z-\lambda|\leq d/4.
\end{equation}
\end{lemma}
\begin{proof}
If $\mathcal{L}_n\subset \Dom (A^2)$, the result has been
established in \linebreak \cite[Lemma~5.3]{bo2}. We treat the more
general case by considering approximate subspace
$\tilde{\mathcal{L}}_n\subset \Dom (A^2)$ with orthonormal bases
sufficiently close to $\mathcal{L}_n$ in the sense specified by
(i)-(iii) below.

As a first step, we recall the following standard result.
For any $v\in \Dom(A)$, there exists a sequence $v_n \in
\Dom(A^2)$ such that $v_n\to v$ and $Av_n \to Av$. That is to say,
$\Dom(A^2)$ is a core (in the operator sense) for $A$.

Let
\begin{displaymath}
c_n = \max\big\{1,\max_{j,k=1...n;\,p,q=0,1}\{\vert\langle
A^pe_j,A^qe_k \rangle\vert\}\big\}.
\end{displaymath}
Then it is always possible to find a set
$\{\tilde{e}_1,\ldots,\tilde{e}_n\}\subset \Dom(A^2)$ such that
\begin{itemize}
\item[(i)] $\{\tilde{e}_1,\ldots,\tilde{e}_n\}$ is orthonormal,
\item[(ii)] $\|e_j - \tilde{e}_j\| \leq c_n^{-1}e^{-n}$,
\item[(iii)]
$|\langle A^pe_j,A^qe_k \rangle - \langle
A^p\tilde{e}_j,A^q\tilde{e}_k \rangle | \leq e^{-n}$,
\end{itemize}
for $j,k=1,\ldots,n$ and $p,q = 0,1$. We may find $\tilde{e}_j$ by
applying the Gram-Schmidt orthogonalisation
procedure to vectors of $\Dom (A^2)$ sufficiently close to the
$e_j$.

Let
$\tilde{\mathcal{L}}_n:=\mathrm{Span}\,\{\tilde{e}_1,\ldots,\tilde{e}_n\}
\subset \Dom (A^2)$. In this proof, the symbol $\sim$ on top of
matrices and operators denotes that they are constructed using
$\mathcal{L}=\tilde{\mathcal{L}}_n$. Note that (iii)
ensures the existence of complex numbers $w^{pq}_{jk}$ such that
$\vert w^{pq}_{jk}\vert\le 1$ and
\begin{displaymath}
\langle A^p\tilde{e}_j,A^q\tilde{e}_k \rangle = \langle
A^pe_j,A^qe_k \rangle + w^{pq}_{jk}e^{-n}.
\end{displaymath}

Property (iii) yields
\begin{align*}
  |[\tilde{M}^{(n)}(z)&-M^{(n)}(z)]_{jk}|  =
  |2z[\tilde{A}_1^{(n)}-A_1^{(n)}]_{jk}+[\tilde{A}_{0}^{(n)}-A_{0}^{(n)}]_{jk}| \\
   & \leq 2|z| |\langle A \tilde{e}_j,\tilde{e}_k \rangle-
    \langle A e_j,e_k \rangle| +
 |\langle A \tilde{e}_j,A \tilde{e}_k \rangle-
    \langle A e_j,A e_k \rangle| \\
  & \leq (2|z|+1)e^{-n},
\end{align*}
Thus,
\begin{equation} \label{e23}
   \|\tilde{M}^{(n)}(z)-M^{(n)}(z)\|\leq (2|z|+1)ne^{-n}.
\end{equation}

Let $u\in \Dom (A)$ be such that $Au=\lambda u$.
We next show that (ii) and the fact that (H) holds for
$\mathcal{L}_n$, ensures that (H) also holds for
$\tilde{\mathcal{L}}_n$. Indeed,
\begin{align*}
\Vert\sum_{k=1}^n\langle & A^p\tilde{\Pi}_n u,A^q\tilde{e}_k\rangle
\tilde{e}_k - \lambda^{p+q} u\Vert \\ & \le \Vert\sum_{k=1}^n \langle
A^p\tilde{\Pi}_n u,A^q\tilde{e}_k\rangle\tilde{e}_k - \langle
A^p\Pi_n
u,A^q e_k\rangle e_k\Vert\\
& \qquad +\Vert\sum_{k=1}^n\langle A^p\Pi_n u,A^q e_k\rangle e_k -
\lambda^{p+q} u\Vert\\
& \le \Vert\sum_{k=1}^n \langle A^p\tilde{\Pi}_n
u,A^q\tilde{e}_k\rangle\tilde{e}_k - \langle A^p\tilde{\Pi}_n
u,A^q\tilde{e}_k\rangle e_k\Vert\\
&\quad +\Vert\sum_{k=1}^n \langle A^p\tilde{\Pi}_n
u,A^q\tilde{e}_k\rangle e_k - \langle A^p\Pi_n u,A^q e_k\rangle
e_k\Vert \\
& \quad +\Vert\sum_{k=1}^n\langle A^p\Pi_n u,A^q e_k\rangle e_k -
\lambda^{p+q} u\Vert
=T_1+T_2+T_3.
\end{align*}
Since $\mathcal{L}_n$ satisfies condition (H), $T_3\to 0$. We must
show that the first two terms also converge to
zero. Consider the first term,
\begin{align*}
T_1 & =\Vert\sum_{jk=1}^n \langle u,\tilde{e}_j\rangle\langle
A^p\tilde{e}_j,A^q\tilde{e}_k\rangle(\tilde{e}_k - e_k)\Vert\\
&\le\Vert u\Vert\sum_{jk=1}^n \vert\langle
A^p\tilde{e}_j,A^q\tilde{e}_k\rangle\vert\Vert\tilde{e}_k - e_k\Vert\\
&=\Vert u\Vert\sum_{jk=1}^n \vert(\langle A^p e_j,A^q
e_k\rangle + w^{pq}_{jk}e^{-n})\vert\Vert\tilde{e}_k - e_k\Vert\\
&\le\Vert u\Vert ne^{-n}\sum_{k=1}^n\Vert\tilde{e}_k -
e_k\Vert+\Vert u\Vert\sum_{jk=1}^n \vert\langle
A^pe_j,A^qe_k\rangle\vert\Vert\tilde{e}_k - e_k\Vert.
\end{align*}
Using (ii) it is clear that $T_1\to0$. For the
second term we have,
\begin{align*}
T_2& =\Vert\sum_{jk=1}^n \langle u,\tilde{e}_j\rangle\langle
A^p\tilde{e}_j,A^q\tilde{e}_k\rangle e_k - \langle
u,e_j\rangle\langle A^p e_j,A^q e_k\rangle e_k\Vert\\
&=\Vert\sum_{jk=1}^n \langle u,\overline{\langle
A^p\tilde{e}_j,A^q\tilde{e}_k\rangle}\tilde{e}_j-\overline{\langle
A^p e_j,A^q
e_k\rangle} e_j\rangle e_k\Vert\\
& \le\Vert u\Vert\sum_{jk=1}^n\Vert\langle
A^p\tilde{e}_j,A^q\tilde{e}_k\rangle\tilde{e}_j-\langle A^p e_j,A^q
e_k\rangle e_j\Vert\\
& =\Vert u\Vert\sum_{jk=1}^n\Vert(\langle A^p e_j,A^q
e_k\rangle + w^{pq}_{jk}e^{-n})\tilde{e}_j-\langle A^p e_j,A^q
e_k\rangle e_j\Vert\\
& \le\Vert u\Vert n^2e^{-n}+\Vert u\Vert\sum_{jk=1}^n \vert\langle
A^pe_j,A^qe_k\rangle\vert\Vert\tilde{e}_j - e_j\Vert.
\end{align*}
Again, using (ii) it is clear that $T_2\to 0$. This ensures that
$\tilde{\mathcal{L}}_n$ also satisfies (H) so
 \eqref{e14} is valid for $\tilde{G}_n(z)$.

The
proof of \eqref{e14} follows. Fix $\varepsilon>0$. Let $\mathbf{v}_n\in
\C^n$ such that $\|\mathbf{v}_n\|=1$ and
\[
   \| M^{(n)}(z) \mathbf{v}_n \| \leq G_n(z)+\varepsilon.
\]
Then, by virtue of \eqref{e23},
\begin{align*}
   \tilde{G}_n(z) &\leq \|\tilde{M}^{(n)} (z) \mathbf{v}_n \| \\
 & \leq \|[\tilde{M}^{(n)}(z)-M^{(n)}(z)]\mathbf{v}_n\| +\| M^{(n)}(z) \mathbf{v}_n \| \\
 & \leq (2|z|+1)ne^{-n} + G_n(z) +\varepsilon.
\end{align*}
As this happens for all $\varepsilon>0$, the fact that
$\tilde{G}_n(z)$ satisfies \eqref{e14} implies that also $G_n(z)$
satisfies this inequality.
\end{proof}

\medskip

\begin{proof}[Proof of Theorem~\ref{t3}]

Let $0<s\leq 1$ be as in Lemma~\ref{t2} and $d$ be as in the hypothesis
of the theorem. By virtue of
Lemmas~\ref{t1} and \ref{t2},
there exists $ N >0$ such that, $G_n(\lambda)\leq
s\delta^2$  and
\begin{equation} \label{e11}
  G_n(z) > s\delta^2\quad
  \mathrm{whenever}\ \delta < |z-\lambda|\leq d/4,
\end{equation}
for all $n>N$. The subharmonicity of $G_n(z)^{-1}$ ensures that the
only local minima of $G_n(z)$ are those points where the function
vanishes. Thus, for all $n>N$, there always exists $\zeta_n\in
\C$ satisfying
\begin{equation} \label{e12}
   |\zeta_n-\lambda|<\delta \qquad \mathrm{and} \qquad
   G_n(\zeta_n)=0.
\end{equation}

Let $\varepsilon_0,\varepsilon_1,\varepsilon_2>0$ be small enough such
that
\[
    \varepsilon_0+2\varepsilon_1|z|+\varepsilon_2|z|^2<s\delta^2,
    \qquad |z-\lambda|<d/4.
\]
Suppose that $n> N$.
Then, by \eqref{e11},
\[
    G_n(z)-(\varepsilon_0+2\varepsilon_1|z|+\varepsilon_2|z|^2)>
    0
\]
for all $\delta< |z-\lambda|\leq d/4$, so
$\Lambda^{\mathcal{L}_n}(\varepsilon_0,\varepsilon_1,\varepsilon_2)\cap
\{\delta < |z-\lambda|\leq d/4 \}=\varnothing$. This, along with
Lemma~\ref{e3}, ensures (b). On the other hand, by
virtue of \eqref{e12},
$\Lambda^{\mathcal{L}_n}(\varepsilon_0,\varepsilon_1,\varepsilon_2)\cap
\{ |z-\lambda|<\delta \}\not =\varnothing$. Thus
Lemma~\ref{t5} yield (a).
\end{proof}


\section{Case Study}
\label{sec4}

Finite rank perturbations of multiplication operators
have been considered previously in connection with spectral pollution
(see \cite{dapl}, \cite{lesh} and \cite{bo2}) due to their
simple structure. In this final section we report on
various numerical experiments we have performed
on a model operator of this type.

Let $e_{k}(x) = (2\pi)^{-\frac{1}{2}}e^{ikx}$ and
\[
   a(x) = \left\{ \begin{array}{ll} 0 & \textrm{for } -\pi\le x<0,\\
   1 & \textrm{for } 0\le x<\pi .
   \end{array} \right.
\]
In this section we assume that
\begin{gather*}
 \mathcal{H}=L^{2}[-\pi
,\pi], \\
\mathcal{L}_{2n+1}= \mathrm{span}\, \{e_{-n}(x),\dots
,e_{n}(x)\} \qquad \mathrm{and} \\
A\phi(x) = a(x)\phi(x) + \langle \phi ,e_{0}\rangle e_{0}(x),
\qquad \phi \in \mathcal{H}.
\end{gather*}
Operator $A$ is bounded and self-adjoint in $\mathcal{H}$.
Moreover, the spectrum $A$ is found explicitly.
Since $A$ is a rank one perturbation of the multiplication operator
by the symbol $a$, Weyl's Theorem ensures that $\Spec_{\ess}(A) =
\mathrm{Range}(a)=\{0,1\}$. On the other hand, the isolated
eigenvalues of finite multiplicity of $A$ are the solutions of
\linebreak
$
\langle(\lambda - a)^{-1}e_{0},e_{0}\rangle \ = 1,
$
\cite{da0}. A straightforward calculation reveals
the two solutions $\lambda^{\pm} = 1\pm \sqrt{2}/2$, which
comprise the discrete spectrum of $A$. The eigenvalue $\lambda^-$
is inside the gap $(0,1)$ of the essential spectrum.

As the symbol $a(x)$ is discontinuous, the Fourier basis  $\{e_k\}$
is not a good choice for approximating $\lambda_-$
using the Galerkin method. Indeed, the solutions
of (P) pollute the whole interval $[0,1]$ as the dimension of
$\mathcal{L}_{2n+1}$ increases. Let us test the quadratic
method described in the preceding sections in this very simple model.

Since $A$ is bounded and $\Pi_{2n+1}\phi\to\phi$
for all $\phi\in \mathcal{H}$, condition (H) of Section~\ref{sec3}
is satisfied.
Thus, by virtue of Theorem~\ref{t3}, both discrete eigenvalues are approached by solutions of (Q) as $n\to\infty$, free from spectral pollution according to Theorem~\ref{pract}.

\smallskip

All the calculation described in this section were carried out using
the computer package MATLAB. Fully functional m-codes
are available at the web page \cite{web}.

We compute the exact solutions of (Q), by finding the $\zeta \in
\C$ such that $\det M^{(2n+1)}(\zeta) = 0$.
The matrix coefficients $A_p^{(2n+1)}$ may be found explicitly using \eqref{e10}. They are sparse and Hermitian with
entries either purely real or purely imaginary. The errors in solving
(Q) are negligible for $n$
of reasonable size  ($<1000$).

In order to test the results established in the
previous sections,
we force large errors in the matrix
entries, and compute the corresponding ``perturbed'' solution of the
problem ($\tilde{\mathrm{Q}}$). For simplicity, we fix
$\varepsilon_0=\varepsilon_1=\varepsilon_2=\varepsilon$.

Let
\begin{equation} \label{e4}
[\tilde{A}_p]_{jk} = [A_p]_{jk} + \frac{\varepsilon}{2n+1}
\alpha_{jk}^{(p)}
\end{equation}
where $\alpha_{jk}^{(p)}$ are random variables sampled from the
unit disk $\{|z|\leq 1\}$ with additional constraints specified below. Then
\begin{align*}
   \|(A_p-\tilde{A}_p)\mathbf{u}\|^2 & =
   \sum_{j=1}^{2n+1}  \left|\sum _{k=1}^{2n+1}
    \frac{\varepsilon}{2n+1}
   \alpha^{(p)}_{jk} u_k \right|^2 \\
    &\leq
   \frac{\varepsilon^2}{(2n+1)^2}
    \sum _{k=1}^{2n+1} |u_k|^2 \sum _{jk=1}^{2n+1}
   |\alpha^{(p)}_{jk}|^2  \\
   & \leq \frac{\varepsilon^2}{(2n+1)^2}
   \sum _{k=1}^{2n+1} |u_k|^2  \sum _{jk=1}^{2n+1} 1
   \leq \varepsilon^2 \sum _{k=1}^{2n+1} |u_k|^2,
\end{align*}
so $\|A_p-\tilde{A}_p\|\leq \varepsilon$. Moreover this bound is
sharp. Indeed, the matrix $T$ such that
$[T]_{jk}=\frac{\varepsilon}{2n+1}$ for all $1\leq j,k \leq 2n+1$,
satisfies
\[
   \|T\|^2 = \|T^\ast T\|= \varepsilon \| T\|.
\]

We consider two types of restrictions on the random variable
$\alpha_{jk}^{(p)}$.
On the one hand, Theorem~\ref{pract} covers the
general situation of moving all entries of $A_p$ along randomly
chosen directions in the complex plane.
Thus, we perform \emph{unstructured perturbations}
by allowing all $\alpha^{(p)}_{jk}\not=0$.
On the other hand, however, in order to reproduce
the effect made by rounding errors in the estimation
of the entries, we perform \emph{non-zero-Hermitian perturbations}
by imposing the condition:
\begin{equation*}
\alpha^{(p)}_{jk}=\left\{\begin{array}{ll}
0 & \mathrm{if\ } [A_p]_{jk}=0,\\
\overline{\alpha^{(p)}_{kj}}\not=0 & \mathrm{if\ }
[A_p]_{jk}\not=0. \end{array} \right.
\end{equation*}

\smallskip

\begin{figure}[p]
\epsfig{file=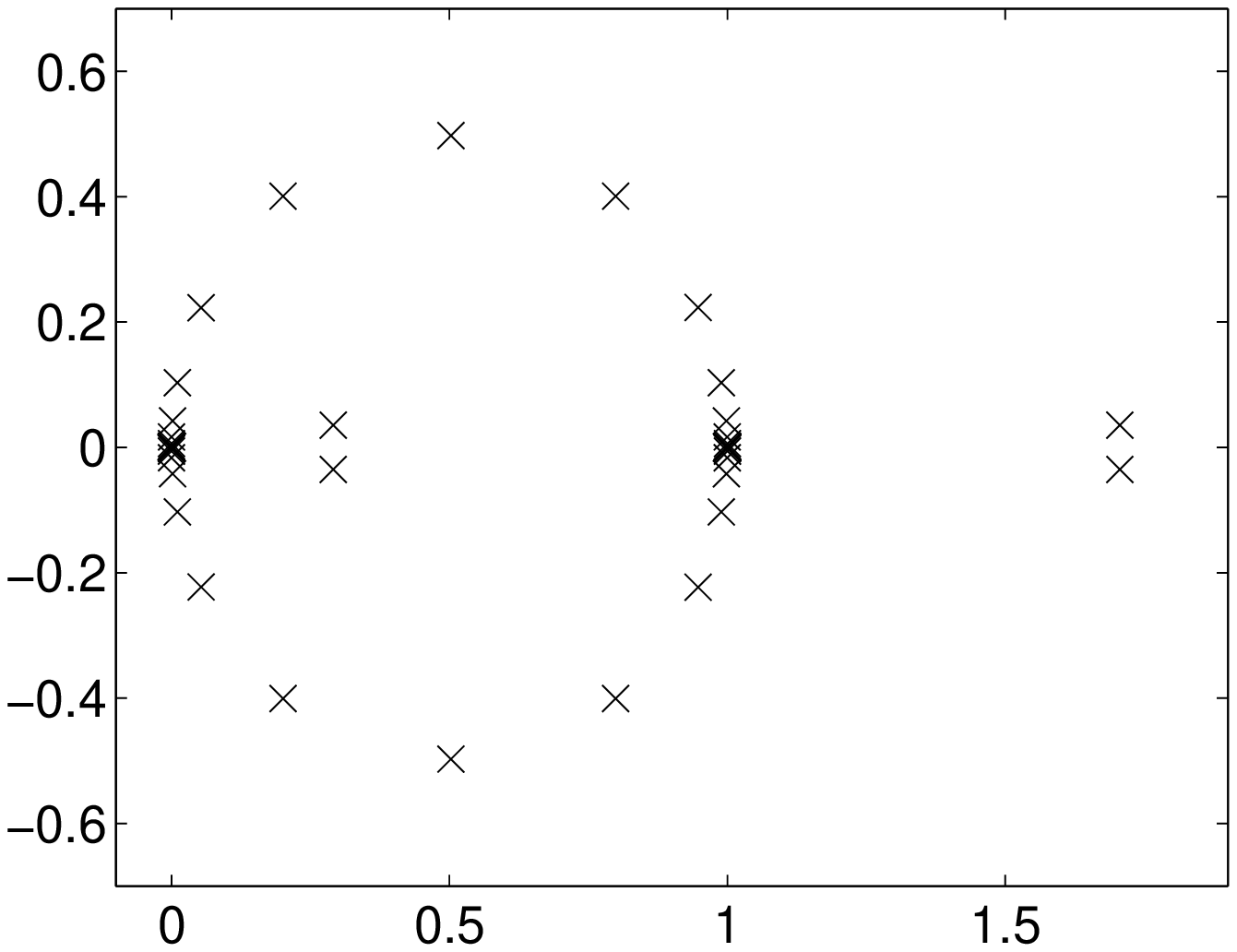, width=9cm}

\caption{Exact solutions to $(\mathrm{Q})$ for $n=50$.}\label{f1}


\epsfig{file=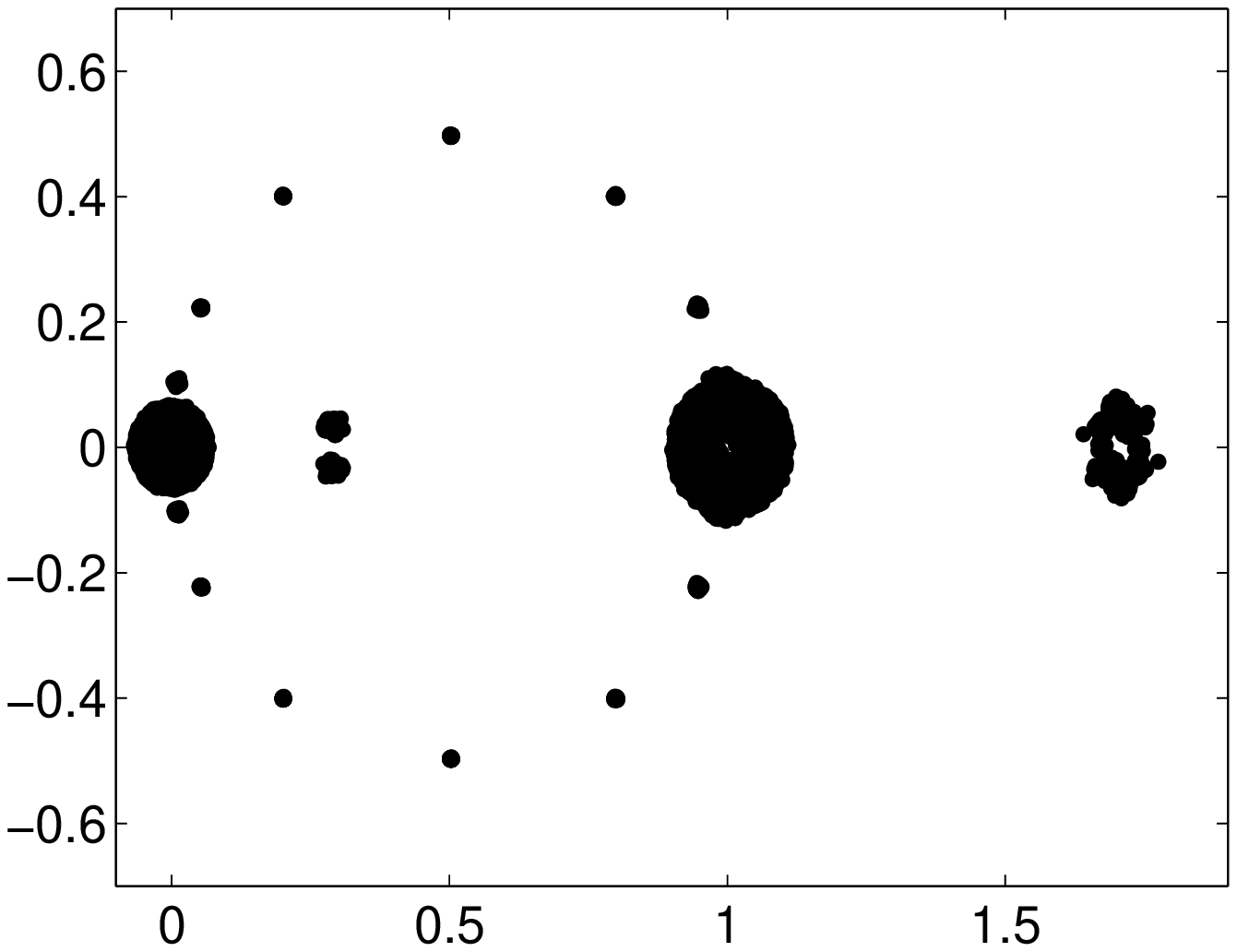, width=9cm}
\epsfig{file=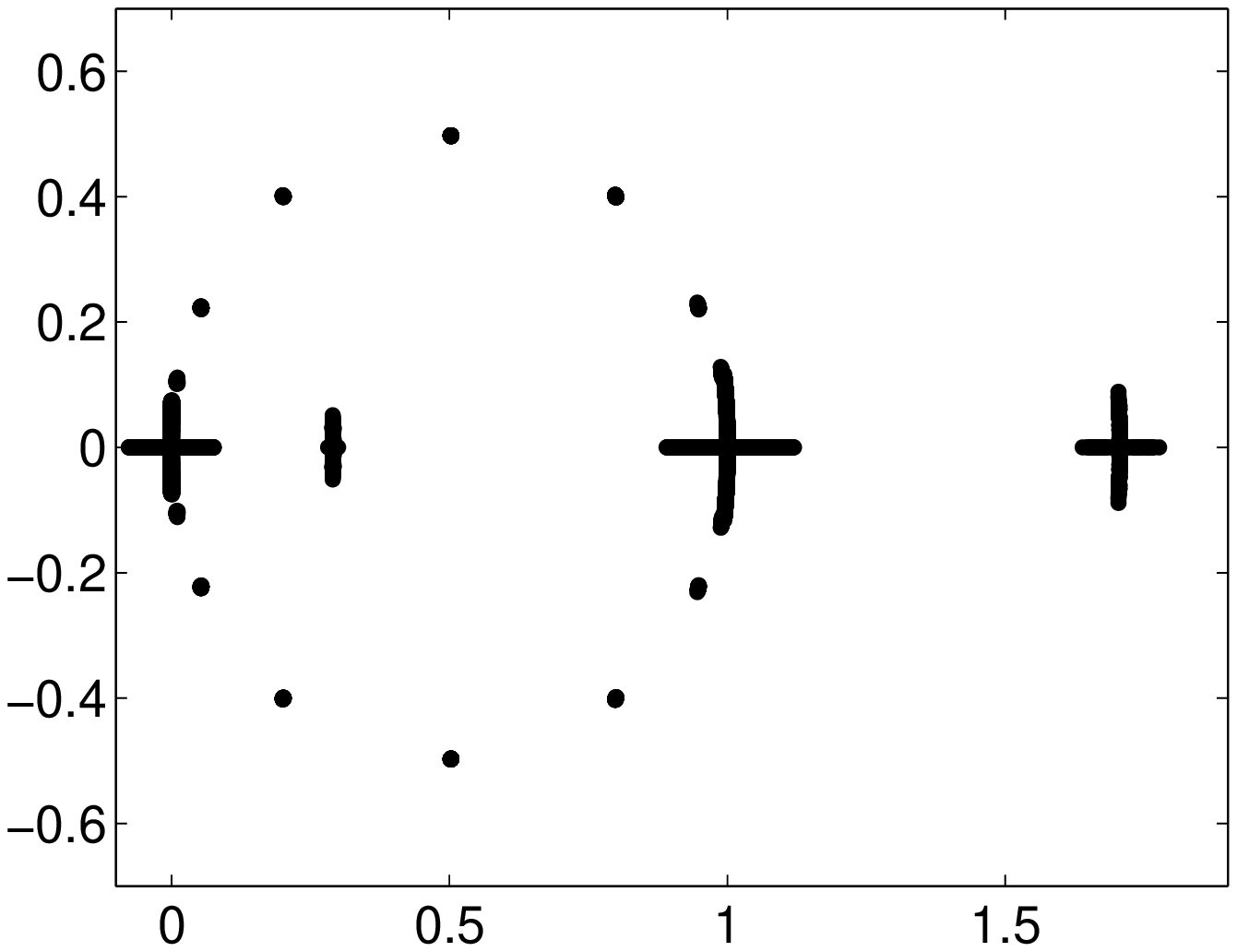, width=9cm}

\caption{Top: solutions to $(\tilde{\mathrm{Q}})$ for 100 unstructured random
perturbations. Bottom:
solutions to $(\tilde{\mathrm{Q}})$ for 100 non-zero-Hermitian
random perturbations. Here $\varepsilon=10^{-1}$} \label{f2}
\end{figure}

In Figure~\ref{f1} we depict the exact solutions of (Q) for
$n=50$. According to \eqref{eu},
the points which are in the vicinity
of the real axis are necessarily close to the spectrum.

Figure~\ref{f2}, on the other hand, depicts the solutions of
($\tilde{\mathrm{Q}}$) corresponding to 100 different random
perturbations. Each of the graphs were constructed by prescribing a
different constraint on the random variables. Here $n=50$ and
$\varepsilon=10^{-1}$. From the perturbed solutions one can identify
$\textrm{Spec}(A)$ less accurately but, once again, without
pollution by virtue of Theorem~\ref{pract}. The
correction $\delta$ of Theorem~\ref{pract}, will
depend on $\zeta$ and $\varepsilon$, but notably not on $n$.
Furthermore, Theorem~\ref{t3} ensures that the clouds observed in
Figure~\ref{f2} will cluster near to each of the exact solutions of
(Q) as $\varepsilon\to 0$.

\smallskip

\begin{figure}[p]

\epsfig{file=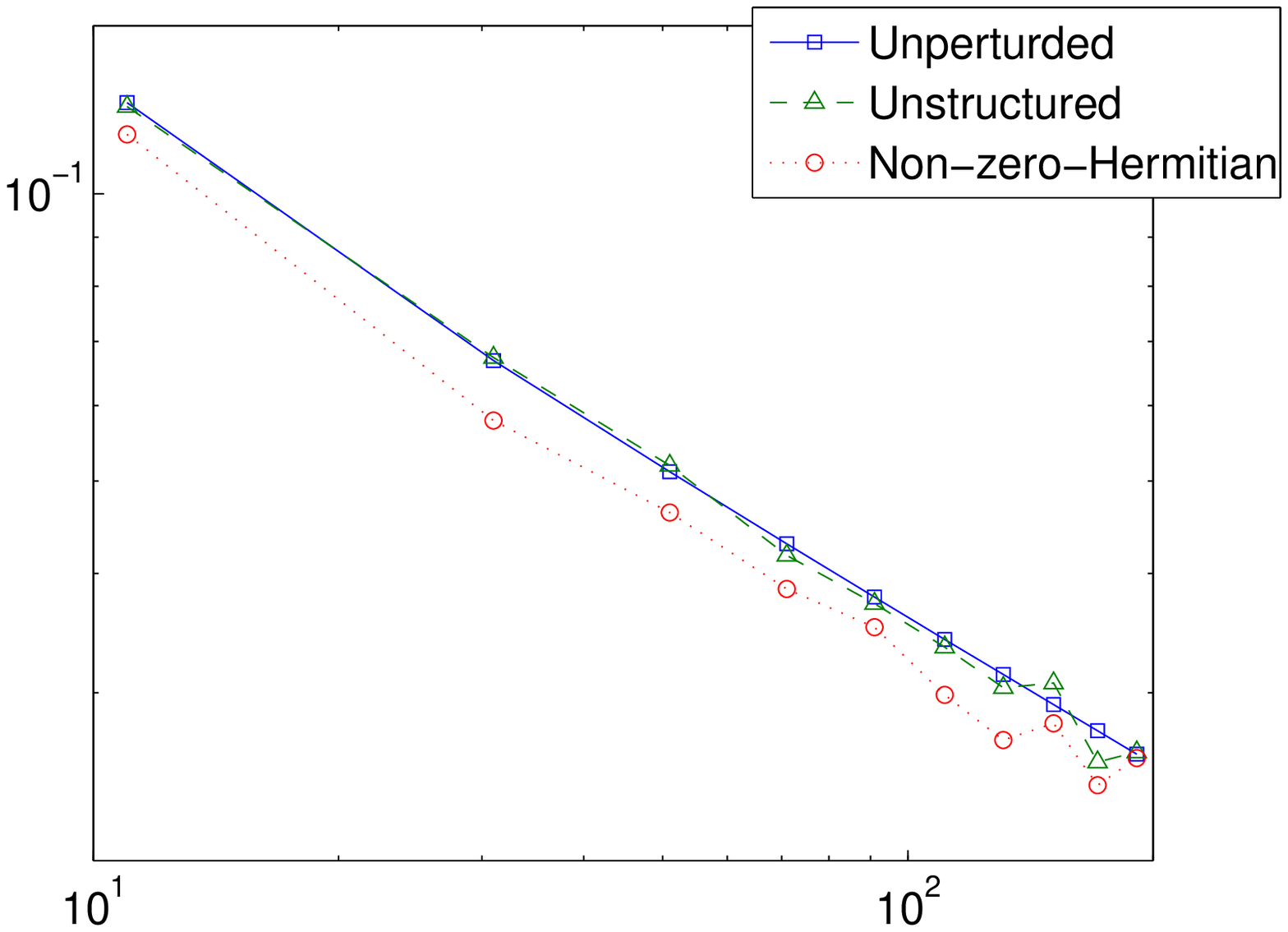, width=9cm}

\caption{Error predicted by Theorem~\ref{pract} in the approximation
of $\lambda^-$.
Here we depict $|\Im \zeta_n^-|$ (unperturbed), $|\Im \zeta_n^{\mathrm{u},-}|$
and $|\Im \zeta_n^{\mathrm{s},-}|$
for $n=5:10:100$. We average the two perturbed
solutions of $(\tilde{\mathrm{Q}})$ over a sample of $20$ problems with
$\varepsilon=10^{-1}$.
The scaling is log-log and the horizontal axis shows $2n+1$.}
\label{f3}

\vspace{1cm}

\epsfig{file=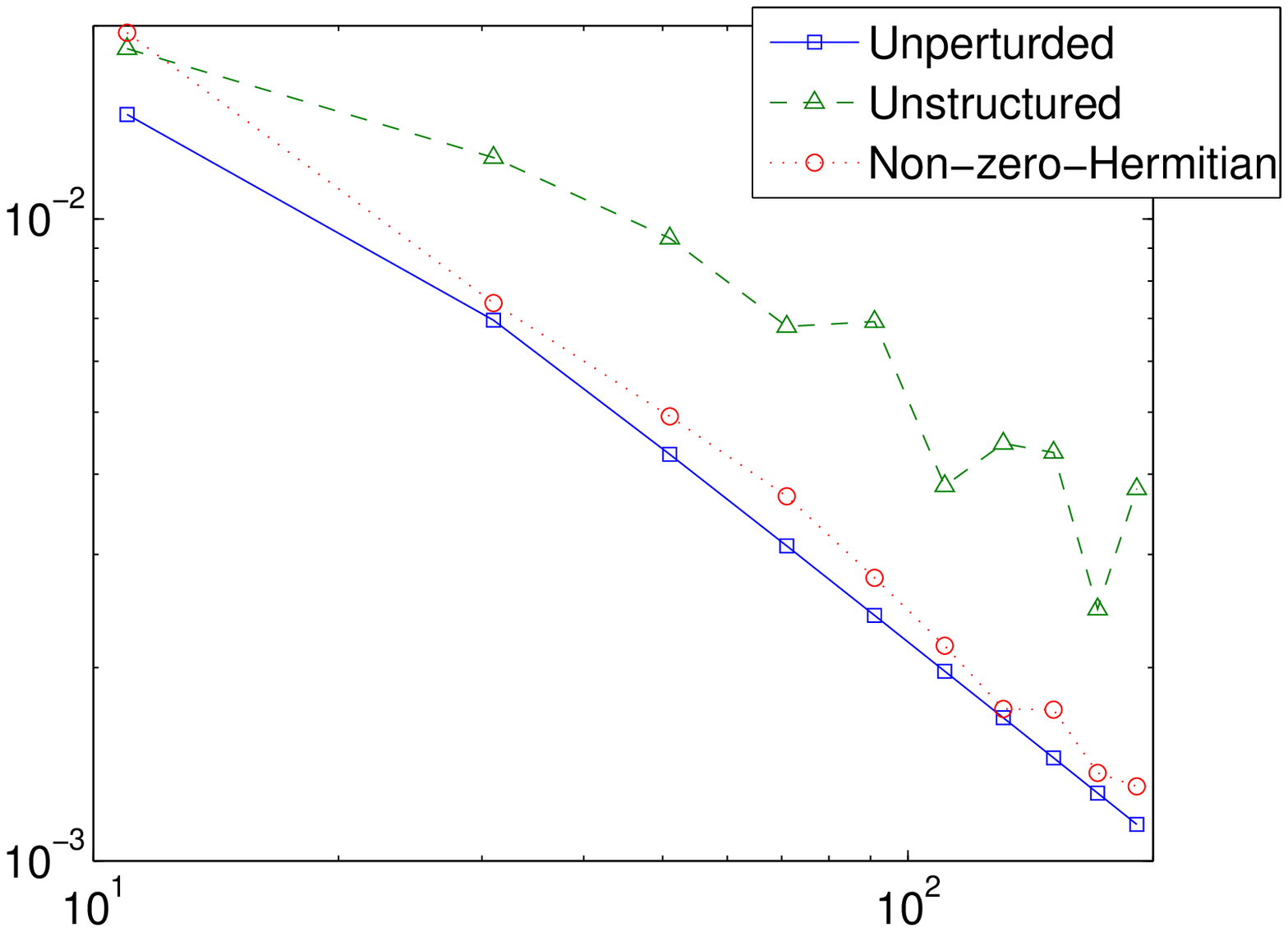, width=9cm}

\caption{Actual error in the approximation
of $\lambda^-$ using the real part of solutions of (Q) and
$(\tilde{\mathrm{Q}})$ .
Here we depict $|\lambda^- - \Re \zeta_n^-|$ (unperturbed), $|\lambda^- - \Re \zeta_n^{\mathrm{u},-}|$
and $|\lambda^- - \Re \zeta_n^{\mathrm{s},-}|$ for $n=5:10:100$. We average the two perturbed solutions of $(\tilde{\mathrm{Q}})$
over a sample of $20$ problems with $\varepsilon=10^{-1}$.
The scaling is log-log and the horizontal axis shows $2n+1$.}
\label{f4}
\end{figure}

Figures~\ref{f3}-\ref{f4} show the outcome of running
Monte Carlo simulations in this model. We fix again
$\varepsilon=10^{-1}$. These pictures have been constructed
in the following manner. For each fixed $n$,
we have found  $\zeta_n^-$, the closest point to the eigenvalue
$\lambda^-$ such that $\det M^{(2n+1)}(\zeta_n^-)=0$.
Then we have performed
$20$ constrained perturbations
and averaged the solutions of $(\tilde{\mathrm{Q}})$
which are closest to $\zeta_n^-$. We know that solutions
of the approximate problems close to $\zeta_n^-$
always exist, as a consequence of Theorem~\ref{t3}. Denote these averages by
$\zeta_n^{\mathrm{u},-}$ and $\zeta_n^{\mathrm{s},-}$ for the unstructured
and non-zero-Hermitian cases respectively. In Figure~\ref{f3} we depict
$|\Im \zeta_n^-|$, $|\Im \zeta_n^{\mathrm{u},-}|$ and $|\Im \zeta_n^{\mathrm{s},-}|$
for $n=5:10:100$. Similarly in Figure~\ref{f4} we depict
$|\lambda^- - \Re \zeta_n^-|$, $|\lambda^- - \Re \zeta_n^{\mathrm{u},-}|$
and $|\lambda^- - \Re \zeta_n^{\mathrm{s},-}|$.

Figure~\ref{f3} provides clear numerical evidence that
the convergence of the quadratic method applied to this simple model
is not lost even when the perturbations are large in modulus.
Figures~\ref{f4} suggests that structured perturbations are
considerably superior to the unstructured ones, in the test
$|\lambda^- - \Re \zeta_n^-|$.

By combining
Figure~\ref{f3} and Theorem~\ref{pract},
we immediately predict a rate
of convergence of $|\lambda^- - \Re \zeta_n^-|=o(n^{-r})$ for $r\approx 1/2$.
It is remarkable, however, that Figure~\ref{f4} strongly suggests
an actual exponent of $r\approx 1$ for this rate of convergence.
An explanation of this phenomenon
is linked to the fact that
$\lambda^-$ is an isolated point of the spectrum, see  \cite[Section~2]{shar}.
We will be reporting on this issue elsewhere.

\section{Acknowledgements} We kindly thank Eugene Shargorodsky
and Michael Levitin for encouraging us to write this manuscript in a first
place and for
their valuable comments during the various stages of its preparation.

\vspace{2cm}

\linespread{1}

\begin{minipage}[c]{15cm}
{\scshape Lyonell Boulton} \\
Department of Mathematics and\\
Maxwell Institute for Mathematical Sciences \\
Heriot-Watt University, Edinburgh EH14 2AS, Scotland\\
E-mail: \texttt{L.Boulton@hw.ac.uk}
\end{minipage}

\vspace{1cm}

\begin{minipage}[c]{15cm}
{\scshape Michael Strauss} \\
Department of Mathematics, Kings College London\\
Strand, London WC2R 2LS, England\\
E-mail: \texttt{michael.strauss@kcl.ac.uk}
\end{minipage}

\newpage

\newpage

\end{document}